\documentclass[12pt,reqno]{amsart}
\usepackage{amsmath,amssymb,amsthm,bbm,fullpage,cite,hyperref,verbatim}

\newtheorem{theorem}{Theorem}[section]
\newtheorem{lemma}[theorem]{Lemma}
\newtheorem{corollary}[theorem]{Corollary}
\newtheorem{observation}[theorem]{Observation}
\theoremstyle{definition}
\newtheorem{definition}[theorem]{Definition}
\newtheorem{remark}[theorem]{Remark}

\newcommand\N{\mathbb{N}}
\newcommand\Z{\mathbb{Z}}
\newcommand\R{\mathbb{R}}
\newcommand\C{\mathbb{C}}
\newcommand\cF{\mathcal{F}}
\newcommand\cG{\mathcal{G}}
\newcommand\cI{\mathcal{I}}
\newcommand\cJ{\mathcal{J}}
\newcommand\id{\mathbbm{1}}
\newcommand\eps{\varepsilon}
\newcommand\ds{\displaystyle}
\newcommand\tg{\tilde s}

\renewcommand\leq{\leqslant}

\renewcommand\le{\leqslant}
\renewcommand\ge{\geqslant}
\newcommand\supp{\operatorname{supp}}

\renewcommand\Im{\operatorname{Im}}
\renewcommand\Re{\operatorname{Re}}

\begin{document}

\title{Flat Littlewood Polynomials Exist}
\author{Paul Balister \and B\'ela Bollob\'as \and Robert Morris \and \\ Julian Sahasrabudhe \and Marius Tiba}

\address{Department of Mathematical Sciences,
University of Memphis, Memphis, TN 38152, USA}
\email{pbalistr@memphis.edu}

\address{Department of Pure Mathematics and Mathematical Statistics, Wilberforce Road,
Cambridge, CB3 0WA, UK, and Department of Mathematical Sciences,
University of Memphis, Memphis, TN 38152, USA}
\email{b.bollobas@dpmms.cam.ac.uk}

\address{IMPA, Estrada Dona Castorina 110, Jardim Bot\^anico, Rio de Janeiro, 22460-320, Brazil}
\email{rob@impa.br}

\address{Peterhouse, Trumpington Street, University of Cambridge, CB2 1RD, UK} 
\email{jdrs2@cam.ac.uk}

\address{Department of Pure Mathematics and Mathematical Statistics, Wilberforce Road,
Cambridge, CB3 0WA, UK}\email{mt576@dpmms.cam.ac.uk}

\thanks{The first two authors were partially supported by NSF grants DMS 1600742 and DMS 1855745,
the third author was partially supported by CNPq (Proc.~303275/2013-8) and FAPERJ (Proc.~201.598/2014),
and the fifth author was supported by a Trinity Hall Research Studentship.}
	
\begin{abstract}
 We show that there exist absolute constants $\Delta > \delta > 0$ such that,
 for all $n \ge 2$, there exists a polynomial $P$ of degree~$n$,
 with coefficients in $\{-1,1\}$, such that
 \[
  \delta\sqrt{n} \le |P(z)| \le \Delta\sqrt{n}
 \]
 for all $z\in\C$ with $|z|=1$. This confirms a conjecture of Littlewood from 1966.
\end{abstract}
	
\maketitle

\section{Introduction}

We say that a polynomial $P(z)$ of degree $n$ is a {\em Littlewood polynomial\/} if 
$$P(z) = \sum_{k=0}^{n} \eps_k z^k,$$
where $\eps_k \in \{-1,1\}$ for all $0 \le k \le n$. The aim of this paper is to prove the following theorem, which answers a question of Erd\H{o}s~\cite[Problem~26]{E57} from 1957, and confirms a conjecture of Littlewood~\cite{L66} from 1966.

\begin{theorem}\label{thm:main}
 There exist constants\/ $\Delta > \delta > 0$ such that, for all\/ $n \ge 2$,
 there exists a Littlewood polynomial\/ $P(z)$ of degree $n$ with 
 \begin{equation}\label{eq:flat:polynomials}
  \delta \sqrt{n} \le |P(z)| \le \Delta \sqrt{n}
 \end{equation}
 for all\/ $z \in \C$ with\/ $|z| = 1$.
\end{theorem}

Polynomials satisfying~\eqref{eq:flat:polynomials} are known as {\em flat polynomials}, and
Theorem~\ref{thm:main} can therefore be more succinctly stated as follows:
``flat Littlewood polynomials exist''. It turns out that our main challenge will be to prove
the lower bound on $|P(z)|$; indeed, explicit polynomials satisfying the upper bound
in~\eqref{eq:flat:polynomials} have been known to exist since the work of Shapiro~\cite{S52}
and Rudin~\cite{R59} over 60 years ago  (see Section~\ref{Shapiro:sec}). 
In the 1980s a completely different (and non-constructive) proof of the upper bound was given by Spencer~\cite{S85}, who used a 
technique that had been developed a few years earlier by Beck~\cite{B81} in his study of combinatorial discrepancy.
We remark that the Rudin--Shapiro construction, and also ideas from discrepancy theory 
(see Section~\ref{brownian:sec}), will play key roles in our proof.

The study of Littlewood polynomials has a long and distinguished history (see, for example,~\cite{B02} or~\cite{M17}), and appears to have originated in the work of Hardy and Littlewood~\cite{HL16} on Diophantine approximation over 100 years ago, in the work of Bloch and P\'olya~\cite{BP32} on the maximum number of real roots of polynomials with restricted coefficients, and in that of Littlewood and Offord~\cite{LO38, LO45, LO48} and others~\cite{EO56,SZ54} on random polynomials. Two important extremal problems that arose from these early investigations are Littlewood's $L_1$-problem~\cite{HL48}, which was famously resolved (up to constant factors) in 1981 by McGehee, Pigno and Smith~\cite{MPS} and Konyagin~\cite{K81}, and Chowla's cosine problem~\cite{C65}, see~\cite{B86,R04}. 

Motivated by this work, Erd\H{o}s~\cite{E57} asked in 1957 whether flat Littlewood polynomials exist, and also, in the other direction, whether there exists a constant $c > 0$ such that, for every polynomial $P_n(z) = \sum_{k=0}^n a_k z^k$ with $a_k \in \C$ and $|a_k| = 1$ for all $ 0 \le k \le n$, we have 
$$|P_n(z)| \ge (1+c) \sqrt{n}$$
for some $z \in \C$ with $|z| = 1$. (Note that, by a simple application of Parseval's theorem, the conclusion holds with $c = 0$.) In the decade that followed, Littlewood wrote a series of papers~\cite{L61,L62a,L62,L64,L66a,L66} on extremal problems concerning polynomials with restricted coefficients. In particular, in~\cite{L66}, and in his book~\cite{L68} on thirty problems in analysis, Littlewood made several conjectures, the best known of which is that flat Littlewood polynomials exist.

\enlargethispage*{\baselineskip}
\enlargethispage*{\baselineskip}
\thispagestyle{empty}

Let us write $\cF_n$ for the family of Littlewood polynomials of degree $n$, and $\cG_n$ for the (larger) family with coefficients satisfying $|a_k| = 1$. The class $\cG_n$ is significantly richer than~$\cF_n$, and for polynomials in this richer class,
significant progress was made in the years following Littlewood's work. It had been known since
the work of Hardy and Littlewood~\cite{HL16} that the upper bound in~\eqref{eq:flat:polynomials}
holds for the polynomial in $\cG_n$ given by setting $a_k := k^{ik}$, and Littlewood~\cite{L62a}
proved that the polynomial in $\cG_n$ given by setting $a_k := \exp\big(\binom{k+1}{2}\pi i/(n+1) \big)$
satisfies the stronger upper and lower bounds
\begin{equation}\label{eq:ultraflat:polynomials}
 |P(z)| = \big( 1 + o(1) \big) \sqrt{n}
\end{equation}
for all $z \in \C$ with $|z| = 1$ except in a small interval around $z = 1$. Following further
progress in~\cite{BN,B77}, and building in particular on work of K\"orner~\cite{Ko80}, the second question
of Erd\H{o}s~\cite{E57} mentioned above was answered by Kahane~\cite{Ka80}, who proved that there exist
\emph{ultra-flat\/} polynomials in~$\cG_n$, i.e., polynomials that satisfy~\eqref{eq:ultraflat:polynomials}
for all $z \in \C$ with $|z| = 1$. More recently, Bombieri and Bourgain~\cite{BB} improved Kahane's bounds,
and moreover gave an effective construction of an ultra-flat polynomial in~$\cG_n$.

For the more restrictive class of Littlewood polynomials, much less progress has been made over the
past 50 years. The Rudin--Shapiro polynomials mentioned above satisfy the upper bound
in~\eqref{eq:flat:polynomials} with $\Delta = \sqrt{2}$ when $n = 2^t-1$, and with $\Delta = \sqrt{6}$
in general (see~\cite{PB}). However, the previously best-known lower bound,
proved by Carroll, Eustice and Figiel~\cite{CEF} via a simple recursive construction, states
that there exist Littlewood polynomials $P_n(z) \in \cF_n$ with $|P_n(z)| \ge n^{0.431}$ for all
sufficiently large $n \in \N$. Moreover, exhaustive search for small values of~$n$ (see~\cite{O18}) suggests
that ultra-flat Littlewood polynomials most likely do not exist. Let us mention one final interesting result in the direction of
Littlewood's conjecture, due to Beck~\cite{B91}, who proved that there exist flat polynomials
in $\cG_n$ with $a_k^{400} = 1$ for every~$k$.

In the next section, we outline the general strategy that we will use to prove Theorem~\ref{thm:main}.
Roughly speaking, our Littlewood polynomial will consist (after multiplication by a suitable negative
power of~$z$) of a real cosine polynomial that is based on the Rudin--Shapiro construction, and an
imaginary sine polynomial that is designed to be large in the (few) places where the cosine polynomial
is small. To be slightly more precise, we will attempt to ``push'' the sine polynomial far away from
zero in these few dangerous places, using techniques from discrepancy theory 
to ensure that we can do so. In order to make this argument work, it will be important that the intervals on which the
Rudin--Shapiro construction is small are ``well-separated'' (see Definition~\ref{def:wellseparated}, below).
The properties of the cosine polynomial that we need are stated in Theorem~\ref{thm:cosine}, and
proved in Section~\ref{Shapiro:sec}; the properties of the sine polynomial are stated in
Theorem~\ref{thm:sine}, and proved in Section~\ref{sine:sec}.


\section{Outline of the proof}\label{outline:sec}

We may assume that $n$ is sufficiently large, since the polynomial $1 - z - z^2 - \dots - z^n$
has no roots with $|z| = 1$ if $n \ge 2$. It will also suffice to prove Theorem~\ref{thm:main}
for $n \equiv 0 \pmod 4$, since the addition of a constant number of terms of the form $\pm z^k$ can
at worst only change $|P(z)|$ by an additive constant. We can also multiply the polynomial by $z^{-2n'}$
so that it becomes the centred `Laurent polynomial'
\[
 \sum_{k = -2n'}^{2n'} \eps_k z^k,
\]
where $n=4n'$. The following theorem therefore implies Theorem~\ref{thm:main}.

\begin{theorem}\label{thm:main:constants}
 For every sufficiently large\/ $n \in \N$, there exists a Littlewood polynomial\/
 $P(z) = \sum_{k = -2n}^{2n} \eps_k z^k$ such that
 \[
  2^{-160} \sqrt{n} \le |P(z)| \le 2^{12} \sqrt{n}
 \]
 for all\/ $z \in \C$ with\/ $|z| = 1$.
\end{theorem}

We remark that the constants in Theorem~\ref{thm:main:constants} could be improved somewhat,
but we have instead chosen to (slightly) simplify the exposition wherever possible.

\subsection{Strategy}\label{strategy:sec}

Before embarking on the technical details of the proof, let us begin by giving a rough outline
of the strategy that we will use to prove Theorem~\ref{thm:main:constants}.
The first idea is to choose a set $C \subseteq [2n] = \{1,\ldots,2n\}$, and set $\eps_{-k} = \eps_k$ for each $k \in C$, and $\eps_{-k} = - \eps_k$ for each $k \in S := [2n] \setminus C$. 
Setting $z = e^{i\theta}$, the polynomial $P(z)$ then decomposes as
\[
 \sum_{k = -2n}^{2n} \eps_k z^k
 = \eps_0 + 2\sum_{k \in C}\eps_k\cos(k\theta)+2i\sum_{k\in S}\eps_k\sin(k\theta).
\]
The real part of this expression is a cosine polynomial, while the imaginary part is a sine polynomial. Our aim is to choose the sine and cosine polynomials so that both are $O(\sqrt{n})$ for all $\theta$, and so that the sine polynomial is large whenever the cosine polynomial is small.  

Let us first describe our rough strategy for choosing the sine polynomial $s(\theta)$, given a suitable cosine polynomial $c(\theta)$. For each `bad' interval $I \subseteq \R / 2\pi\Z$ on which $|c(\theta)| < \delta \sqrt{n}$, we will choose a direction (positive or negative), and attempt to `push' the sine polynomial in that direction on that interval. In other words, we pick a step function that is $\pm K \sqrt{n}$ on each of the bad intervals, and zero elsewhere, where $K$ is a large constant. We then attempt to approximate this step function with a sine polynomial, the hope being that we can do so with an error of size $O(\sqrt{n})$ on each bad interval (independent of $K$).

In order to carry out this plan, we will use an old result\footnote{We will in fact find it convenient to use a variant of Spencer's theorem, due to Lovett and Meka~\cite{LM}.} of Spencer~\cite{S85} on combinatorial discrepancy (in the form of Corollary~\ref{cor:colouring} below), first to choose the step function, and then to show that we can approximate it sufficiently closely. More precisely, the first application (see Lemma~\ref{lem:choosing:alpha}) provides us with a step function whose Fourier coefficients are all small, and the second application (see Lemmas~\ref{lem:sin_derivatives} and~\ref{lem:choosing:epsilon}) then produces a sine polynomial that does not deviate by more that $O(\sqrt{n})$ from this step function. 

To make the sketch above rigorous, we will need the bad intervals to have a number of useful properties; roughly speaking, they should be `few', `small', and `well-separated'. In particular, we will construct (see Definition~\ref{def:wellseparated} and Theorem~\ref{thm:cosine}) a set $\cI$ of intervals, each of size $O(1/n)$, separated by gaps of size $\Omega(1/n)$, with $|c(\theta)| \ge \delta \sqrt{n}$ for all $\theta \notin \bigcup_{I \in \cI} I$. Moreover, the number of intervals in $\cI$ will be at most $\gamma n$ for some small constant $\gamma > 0$. 

To see that these demands are not unreasonable, note first that if $C \subseteq [\gamma n]$ then the cosine polynomial has few roots, and the `typical' value of the derivative of the cosine polynomial should be $\Theta((\gamma n)^{3/2})$. This means that, if we choose $\delta$ much smaller than $\gamma$, the polynomial should typically vary by more than $\delta \sqrt{n}$ over a distance of order $1/n$. In particular, we will show that if the set of bad intervals cannot be covered by a collection of small and well-separated intervals (in the sense described above), then several of the derivatives must be small simultaneously. For our cosine polynomial we shall use an explicit construction based on the Rudin--Shapiro polynomials (see Section~\ref{Shapiro:sec}), and we will show (see Lemma~\ref{lem:H_deriv_lbounds}) that the value and first three derivatives of this polynomial cannot all be simultaneously small. 

\enlargethispage*{\baselineskip}
\thispagestyle{empty}

\subsection{The cosine polynomial}

Let $n \in \N$ be sufficiently large, choose $2^{-43} < \gamma \le 2^{-40}$ such that
\begin{equation}\label{eq:gamma}
 \gamma n = 2^{t+11} + 2^t - 1
\end{equation}
for some {\em odd\/} integer~$t$, and set
\[
 \delta := 2^{-8} \gamma^{7/2},
\]
noting that $\delta>2^{-160}$. Define $C \subseteq [2\gamma n]$ by setting $C = 2C'$, where
\[
 C' := \big\{ 2^{t+10}, \dots, 2^{t+10} + 2^t - 1 \big\} \cup \big\{ 2^{t+11}, \dots, 2^{t+11} + 2^t - 1 \big\},
\]
so that $C$ is a set of $2^{t+1}$ {\em even\/} integers. Our first aim is to construct a cosine polynomial
\[
 c(\theta) = \sum_{k \in C} \eps_k \cos(k\theta),
\]
that is only small on a few, well-separated intervals, and is never too large. 

To state the two main steps in the proof of Theorem~\ref{thm:main:constants}, we first need to define what we mean by a `suitable' and `well-separated' family of intervals.

\begin{definition}\label{def:wellseparated}
 Let $\cI$ be a collection of disjoint intervals in $\R/2\pi\Z$. We will say that $\cI$ is \emph{suitable\/} if
 \begin{itemize}
  \item[$(a)$] The endpoints of each interval in $\cI$ lie in $\frac{\pi}{n} \Z$;
  \item[$(b)$] $\cI$ is invariant under the maps $\theta \mapsto \pi \pm \theta$;
  \item[$(c)$] $|\cI| = 4N$ for some $N \le \gamma n$.
 \end{itemize}
 We say that a suitable collection $\cI$ is \emph{well-separated\/} if
 \begin{itemize}
  \item[$(d)$] $|I| \le 6\pi/n$ for each $I \in \cI$;
  \item[$(e)$] $d(I,J) \ge \pi/n$ for each $I,J \in \cI$ with $I \ne J$;\footnote{Given two sets
   $I,J \subseteq \R/2\pi\Z$, let us write $d(I,J) := \inf\{ d(\theta,\theta') : \theta \in I,\, \theta'\in J\}$, where $d(\theta,\theta')$ is the distance between $\theta$ and $\theta'$ mod $2\pi$.}
  \item[$(f)$] $\bigcup_{I \in \cI}I$ is disjoint from the set $(\pi/2)\Z + [-100\pi/n,100\pi/n]$.
 \end{itemize}
\end{definition}

We will prove the following theorem about cosine polynomials.

\begin{theorem}\label{thm:cosine}
 There exists a cosine polynomial
 \[
  c(\theta) = \sum_{k \in C} \eps_k \cos(k\theta),
 \]
 with\/ $\eps_k \in \{-1,1\}$ for every\/ $k \in C$, and a suitable and well-separated collection\/
 $\cI$ of disjoint intervals in\/ $\R/2\pi\Z$, such that
 \[
  |c(\theta)| \ge \delta \sqrt{n}
 \]
 for all\/ $\theta \notin \bigcup_{I \in \cI}I$, and\/ $|c(\theta)|\le \sqrt{n}$
 for all\/ $\theta \in \R/2\pi\Z$.
\end{theorem}

The cosine polynomial we will use to prove Theorem~\ref{thm:cosine} is a slight modification
of the Rudin--Shapiro polynomial.  We might remark here that one would expect almost any cosine
polynomial whose absolute value is $O(\sqrt{n})$ to satisfy somewhat similar conditions,
but this seems difficult to prove in general.

\subsection{The sine polynomials}

There will in fact be two sine polynomials; the first,
\begin{equation}\label{def:s:even:sketch}
 s_e(\theta)=\sum_{j\in S_e}\eps_j\sin(j\theta),
\end{equation}
will just be chosen to be small everywhere, more precisely at most $6\sqrt{n}$ for all $|z| = 1$
(see Lemma~\ref{lem:even_ubounds} below). It is defined on the set $S_e=2S'_e$, where
\[
 S'_e := [n] \setminus C'
\]
so that $S_e$ is the set of remaining even integers in $[2n]$.

We write $S_o := \{1,3,\dots,2n-1\}$ for the set of all the odd integers in $[2n]$, and our main
task will be to construct an `odd sine polynomial'
\[
 s_o(\theta) = \sum_{k \in S_o} \eps_k \sin(k\theta)
\]
that is large on each $I \in \cI$, and not too large elsewhere.
To be precise, we shall prove the following theorem.

\begin{theorem}\label{thm:sine}
 Let\/ $\cI$ be a suitable and well-separated collection of disjoint intervals in\/ $\R/2\pi\Z$.
 There exists a sine polynomial
 \[
  s_o(\theta) = \sum_{k \in S_o} \eps_k \sin(k\theta),
 \]
 with\/ $\eps_k \in \{-1,1\}$ for every\/ $k \in S_o$, such that
 \begin{itemize}
  \item[$(i)$] $|s_o(\theta)|\ge 10\sqrt{n}$ for all\/ $\theta\in\bigcup_{I\in\cI}I$, and
  \item[$(ii)$] $|s_o(\theta)|\le 2^{10} \sqrt{n}$ for all\/ $\theta\in\R$.
 \end{itemize}
\end{theorem}

To deduce Theorem~\ref{thm:main:constants} from the results above, we simply set
\[
 P( e^{i\theta} ) := \big( 1 + 2c(\theta) \big) + 2i \big( s_e(\theta) + s_o(\theta) \big),
\]
where $c(\theta)$ and $s_o(\theta)$ are the cosine and sine polynomials given by Theorems~\ref{thm:cosine} and~\ref{thm:sine} respectively, and $s_e(\theta)$ is a sine polynomial as in~\eqref{def:s:even:sketch} (see Section~\ref{sine:sec} for the details).

The rest of the paper is organised as follows. First, in Section~\ref{Shapiro:sec}, we will define $c(\theta)$
and $s_e(\theta)$, and prove Theorem~\ref{thm:cosine}. In Section~\ref{brownian:sec} we will recall
the main lemma from~\cite{LM} and deduce Corollary~\ref{cor:colouring}; this will be our main tool in the proof of Theorem~\ref{thm:sine}, which is given in Section~\ref{sine:sec}. Finally, we will conclude by completing the proof of Theorem~\ref{thm:main:constants}.

\section{Rudin--Shapiro Polynomials}\label{Shapiro:sec}

In this section we will define the cosine polynomial that we will use to prove
Theorem~\ref{thm:cosine}, and the sine polynomial that we will use on the remaining even
integers. In both cases, we use the so-called Rudin--Shapiro polynomials, which were
introduced independently by Shapiro~\cite{S52} and Rudin~\cite{R59} (and whose sequence
of coefficients was also previously studied by Golay~\cite{G49}). These polynomials
have been extensively studied over the last few decades, see,
e.g.,~\cite{B73,BC,BLM,R17}. Let us begin by recalling their definition.

\begin{definition}[Rudin--Shapiro polynomials]\label{def:Shapiro}
 Set $P_0(z)=Q_0(z)=1$ and inductively define
 \begin{align*}
  P_{t+1}(z)&=P_t(z) + z^{2^t}Q_t(z), \text{ and}\\
  Q_{t+1}(z)&=P_t(z) - z^{2^t}Q_t(z),
 \end{align*}
 for each $t \ge 0$.
\end{definition}

Observe that $P_t(z)$ and $Q_t(z)$ are both Littlewood polynomials of
degree $2^t-1$. A simple induction argument (see, e.g.,~\cite{M17})
shows that $P_t(z)P_t(1/z)+Q_t(z)Q_t(1/z)=2^{t+1}$
for all $z\in\C\setminus\{0\}$. It follows that
\begin{equation}\label{e:Shapiro_norms}
 |P_t(z)|^2 + |Q_t(z)|^2 = 2^{t+1},
\end{equation}
and hence $|P_t(z)|,|Q_t(z)| \le 2^{(t+1)/2}$, for every $z \in \C$ with $|z|=1$.
Observing that the first $2^t$ terms of $P_{t+1}$ are the same as for~$P_t$, let us write
$P_{<n}(z)$ for the polynomial of degree $n-1$ that agrees with $P_t(z)$ on the first $n$
terms for all sufficiently large~$t$, and note that $P_t(z) = P_{<2^t}(z)$. The following bound,
which is a straightforward consequence of~\eqref{e:Shapiro_norms}, was proved by Shapiro~\cite{S52}.
(Stronger bounds are known, see~\cite{PB}, but we shall not need them.)

\begin{lemma}\label{lem:general_Shapiro_bound}
 $|P_{<n}(z)| \le 5\sqrt{n}$ for every\/ $z \in \C$ with\/ $|z| = 1$.
\end{lemma}

We now set
\[
  T := 2^{t + 10},
\]
and define our cosine polynomial to be
\begin{equation}\label{def:c}
 c(\theta) := \Re\big( z^{T} P_t(z) + z^{2T} Q_t(z) \big),
\end{equation}
and our even sine polynomial to be
\begin{equation}\label{def:se}
 s_e(\theta) := \Im\big( P_{<(n+1)}(z) - z^{T} P_t(z) - z^{2T} P_t(z) \big),
\end{equation}
where in both cases $z = e^{2i\theta}$ (note the factor of 2 in the exponent here).
We claim first that\footnote{We define the \emph{support}, $\supp(f)$, of a sine polynomial
$f(\theta)=\sum_{k>0}\eps_k\sin(k\theta)$ or cosine polynomial
$f(\theta)=\sum_{k>0}\eps_k\cos(k\theta)$ to be the set of $k$ such that $\eps_k\ne0$.}
$\supp(c) = C$ and $\supp(s_e) = S_e$.
This is clear for $c$, since $C = 2C'$ and $C' = \{T,\dots,T+2^t-1\} \cup \{ 2T,\dots,2T+2^t-1\}$;
for $s_e$ it follows since the terms of $P_{<(n+1)}(z)$ corresponding to $C'$ form the polynomial
$z^{T} P_t(z) + z^{2T} P_t(z)$. (To see this, simply consider the first time that these terms
appear in Definition~\ref{def:Shapiro}, and note that the first $2^t$ terms of both $P_{t+10}$ and
$Q_{t+10}$ are the same as for~$P_t$.) We remark that the idea behind the definition of $c(\theta)$ is that the
highly oscillatory factors $z^{T}$ and $z^{2T}$ allow us to show that $c$ and its first three derivatives cannot all simultaneously be small (see Lemma~\ref{lem:H_deriv_lbounds}, below).

The following lemma is an almost immediate consequence of Lemma~\ref{lem:general_Shapiro_bound}
and~\eqref{e:Shapiro_norms}.

\begin{lemma}\label{lem:even_ubounds}
 $|c(\theta)| \le \sqrt{n}$ and\/ $|s_e(\theta)| \le 6\sqrt{n}$ for every\/ $\theta\in\R$.
\end{lemma}
\begin{proof}
Observe first that, setting $z := e^{2i\theta}$, we have
\[
 |c(\theta)| \le |P_t(z)| + |Q_t(z)| \le 2^{(t+3)/2} \le \sqrt{n},
\]
where the first inequality follows from the definition of~$c$, the second holds
by~\eqref{e:Shapiro_norms}, and the last holds by~\eqref{eq:gamma}, since $\gamma \le 1$.
Similarly, we have
\[
 |s_e(\theta)| \le |P_{<(n+1)}(z)| + 2|P_t(z)| \le 5\sqrt{n+1} + 2^{(t+3)/2} \le 6\sqrt{n},
\]
where the first inequality follows from the definition of~$s_e$, the second holds by
Lemma~\ref{lem:general_Shapiro_bound} and~\eqref{e:Shapiro_norms}, and the last holds by~\eqref{eq:gamma}.
\end{proof}

In order to prove Theorem~\ref{thm:cosine}, it remains to show that $|c(\theta)| \ge \delta \sqrt{n}$
for all $\theta \notin \bigcup_{I \in \cI} I$, for some suitable and well-separated collection $\cI$
of disjoint intervals in $\R/2\pi\Z$. When doing so we will find it convenient to rescale the
polynomial as follows: define a function $H \colon \R \to \C$ by setting
\[
 H(x) := e^{ix} \alpha(x) + e^{2ix}\beta(x),
\]
where
\[
 \alpha(x) := 2^{-(t+1)/2} P_t( e^{ix/T}) \qquad \text{and} \qquad
 \beta(x)  := 2^{-(t+1)/2} Q_t( e^{ix/T} ),
\]
and observe that
\[
 c(\theta)=2^{(t+1)/2}\,\Re\big( H(2T\theta) \big).
\]
Note that, by~\eqref{e:Shapiro_norms}, we have
\[
 |\alpha(x)|^2 + |\beta(x)|^2 = 1.
\]
We think of $\alpha(x)$ and $\beta(x)$ as being slowly varying functions, relative to the much more
rapidly varying exponential factors in the definition of $H(x)$.

The key property of the polynomial $c(\theta)$ that we will need is given by the following lemma.

\begin{lemma}\label{lem:intervals:nice}
Let\/ $0 < \eta < 2^{-11}$. Every interval\/ $I \subseteq \R$ of length\/~$7\eta$ contains a sub-interval\/ $J \subseteq I$ of length\/ $\eta$ such that
 \[
  \big|\Re\big( H(x) \big) \big| \ge \frac{\eta^3}{2^7}
 \]
 for every\/ $x \in J$. Moreover, if $I=[a,a+7\eta]$ then we can take $J = [a+j\eta,a+(j+1)\eta]$ for some\/ $j\in\{0,1,\dots,6\}$.
\end{lemma}

To prove Lemma~\ref{lem:intervals:nice}, we will first need to prove the following lemma.

\begin{lemma}\label{lem:H_deriv_lbounds}
 For any\/ $x\in\R$ there exists\/ $k \in \{0,1,2,3\}$ such that
 \[
  \big| \Re\big( H^{(k)}(x) \big) \big| \ge \frac{1}{4}.
 \]
\end{lemma}

The proof of Lemma~\ref{lem:H_deriv_lbounds} is not very difficult, but we will need to work a little. We will use Bernstein's classical inequality (see, e.g.,~\cite{S41}), which states that if $f(z)$ is a polynomial of degree $n$, then 
\begin{equation}\label{eq:Bernstein}
\max_{|z| = 1} |f'(z)| \leq n \cdot \max_{|z| = 1} |f(z)|.
\end{equation}
This easily implies the following bound on the derivatives of the Rudin--Shapiro polynomials.  

\begin{lemma}\label{lem:Shapiro_bounds}
 Let\/ $0 \le k, t \in \Z$. We have
\begin{equation}\label{eq:Shapiro:derivatives}
  \bigg| \frac{d^k}{d\theta^k} P_t(e^{i\theta}) \bigg|,\,\bigg|\frac{d^k}{d\theta^k}Q_t(e^{i\theta})\bigg|
  \le \/ 2^{kt + (t+1)/2}
\end{equation}
 for every\/ $\theta\in\R$. In particular, 
\begin{equation}\label{eq:alpha:beta:derivatives}
 |\alpha^{(k)}(x)|,\, |\beta^{(k)}(x)| \le \/ 2^{-10k}
\end{equation}
for every $k \ge 1$ and $x \in \R$.
\end{lemma}

Note that~\eqref{eq:alpha:beta:derivatives} justifies our intuition that $\alpha(x)$ and $\beta(x)$
vary relatively slowly. 

\begin{proof}
To prove~\eqref{eq:Shapiro:derivatives} we simply apply~\eqref{eq:Bernstein} $k$ times, and~\eqref{e:Shapiro_norms} once. It follows from~\eqref{eq:Shapiro:derivatives} that
$$ \max\big\{ |\alpha^{(k)}(x)|,\, |\beta^{(k)}(x)| \big\} \le 2^{-(t+1)/2}
 \cdot T^{-k} \cdot 2^{kt +(t+1)/2} = 2^{-10k}$$
 for every $k \ge 1$ and $x \in \R$, as claimed.
\end{proof}

We will use the following easy consequences of Lemma~\ref{lem:Shapiro_bounds}.

\begin{lemma}\label{lem:H_deriv_ubounds}
 For each\/ $0 \le k \le 4$, and every\/ $x \in \R$, we have
 \[
  \big|H^{(k)}(x)-\big(i^ke^{ix}\alpha(x)+(2i)^ke^{2ix}\beta(x)\big)\big| \le \frac{1}{8}
 \]
 and
 \[
  |H^{(k)}(x)|\le 2^k+2.
 \]
\end{lemma}

\begin{proof}
Since $H(x) = e^{ix} \alpha(x) + e^{2ix}\beta(x)$, we have
\[
 H^{(k)}(x) = \sum_{j=0}^k\binom{k}{j}\big(i^{k-j} e^{ix} \alpha^{(j)}(x) + (2i)^{k-j} e^{2ix} \beta^{(j)}(x)\big),
\]
and hence, using~\eqref{eq:alpha:beta:derivatives},
\[
 \big|H^{(k)}(x) - \big( i^k e^{ix} \alpha(x) + (2i)^k e^{2ix} \beta(x) \big) \big|
 \le \sum_{j=1}^k \binom{k}{j} \big(1 + 2^{k-j} \big) 2^{-10j} \le \frac{1}{8}
\]
(with room to spare) since $k \le 4$.
Since $|i^ke^{ix}\alpha(x)+(2i)^ke^{2ix}\beta(x)|\le 1+2^k$,
it follows immediately that
\[
 |H^{(k)}(x)|\le 2^k+2,
\]
as claimed.
\end{proof}

We can now easily deduce Lemma~\ref{lem:H_deriv_lbounds}.

\begin{proof}[Proof of Lemma~\ref{lem:H_deriv_lbounds}]
Suppose that
\[
 \big| \Re\big( H^{(k)}(x) \big) \big| < \frac{1}{4}
\]
for each $k \in \{0,1,2,3\}$. Setting
\[
 E_k := \Re\big( i^k e^{ix} \alpha(x) + (2i)^k e^{2ix} \beta(x) \big),
\]
observe that
\[\begin{array}{lll}
 \Re\big( e^{ix}\alpha(x) \big) = \ds\frac{4E_0 + E_2}{3},
 && \Re\big( e^{2ix} \beta(x) \big) = - \ds\frac{E_0 + E_2}{3},\\[+2ex]
 \Im\big( e^{ix}\alpha(x) \big) = - \ds\frac{4E_1 + E_3}{3},
 & \qquad \text{and} \qquad
 &\Im\big( e^{2ix}\beta(x) \big) = \ds\frac{E_1 + E_3}{6}.
\end{array}\]
Now, by Lemma~\ref{lem:H_deriv_ubounds}, we have
\[
 |E_k| \le \frac{1}{4} + \frac{1}{8} = \frac{3}{8}
\]
for each $k \in \{0,1,2,3\}$, and therefore
\begin{align*}
 1 & \, = \, |\alpha(x)|^2+|\beta(x)|^2 \\
 & \, = \, \big| \Re\big( e^{ix} \alpha(x) \big)\big|^2 + \big|\Im\big( e^{ix}\alpha(x) \big) \big|^2
 + \big| \Re\big( e^{2ix}\beta(x) \big) \big|^2 + \big|\Im\big( e^{2ix} \beta(x) \big) \big|^2\\
 & \, \le \, \bigg( \frac{5^2}{3^2} + \frac{5^2}{3^2} + \frac{2^2}{3^2} + \frac{2^2}{6^2} \bigg) \cdot \frac{3^2}{8^2}
 \, = \, \frac{55}{9} \cdot \frac{9}{64} \, < \, 1,
\end{align*}
which is a contradiction. It follows that $|\Re( H^{(k)}(x_0))| \ge 1/4$ for some $0 \le k \le 3$.
\end{proof}

To deduce Lemma~\ref{lem:intervals:nice} from Lemmas~\ref{lem:H_deriv_lbounds} and~\ref{lem:H_deriv_ubounds},
we shall use a generalization of Lagrange interpolation from~\cite[Theorem 2]{H91} that
bounds the higher derivatives of a function in terms of its values at certain points.

\begin{theorem}\label{thm:Lagrange}
 Let\/ $f\colon I\to\R$ be a\/ $k+1$ times continuously differentiable function
 and suppose\/ $y_0,\dots,y_k\in I$ with\/ $y_0 < y_1 < \dots < y_k$. Then\footnote{In the notation
 of \cite{H91}, the sum in the first $\|\cdot\|_\infty$ expression is $L^{(k)}(x)$ where
 $L(x)=\sum_i f(y_i) \prod_{j\ne i}(x-y_j)/(y_i-y_j)$, and the second $\|\cdot\|_\infty$ expression is
 $\|\omega^{(k)}(x)/(k+1)!\|_\infty$, where $\omega(x)=\prod_i (x-y_i)$. Note that the inequality is tight
 when $f(x)=\omega(x)$.}
 \[
  \bigg\| f^{(k)}(x) - \sum_{i=0}^k \frac{k! f(y_i)}{\prod_{j\ne i}(y_i - y_j)} \bigg\|_\infty
  \le \ \bigg\| x - \frac{1}{k+1}\sum_{i=0}^k y_i\bigg\|_\infty \cdot \|f^{(k+1)}(x)\|_\infty.
 \]
\end{theorem}

Lemma~\ref{lem:intervals:nice} is a straightforward consequence of Lemmas~\ref{lem:H_deriv_lbounds} and~\ref{lem:H_deriv_ubounds} and Theorem~\ref{thm:Lagrange}. 

\begin{proof}[Proof of Lemma~\ref{lem:intervals:nice}]
Let $I = [a,a+7\eta]$, and suppose (for a contradiction) that for each $0 \le j \le 6$,
there exists a point
\[
 x_j \in I_j := \big[ a + j\eta, \, a + (j+1)\eta \big]
\]
such that $|\Re(H(x_j))| < 2^{-7} \eta^3$. We will show that $|\Re(H^{(k)}(x_0))| < 1/4$
for each $0 \le k \le 3$, which will contradict Lemma~\ref{lem:H_deriv_lbounds}, and hence prove the lemma.

For $k = 0$, we have $|\Re(H^{(k)}(x_0))|<2^{-7}\eta^3<1/4$ (by assumption), so let $k\in\{1,2,3\}$.
By Lemma~\ref{lem:H_deriv_ubounds} and Theorem~\ref{thm:Lagrange}, applied with $f := \Re(H)$
and $y_j := x_{2j}$ for each $0 \le j \le k$ (so, in particular, $|y_i - y_j| \ge \eta$ for all $i\ne j$),
we have
\begin{align*}
 \big|\Re\big( H^{(k)}(x_0) \big)\big|
 & \,\le\, \sum_{i=0}^k\frac{k!}{\eta^k}\cdot \frac{\eta^3}{2^7}
 +7\eta\cdot\big\|\Re\big(H^{(k+1)}(x)\big)\big\|_\infty\\
 & \,\le\, \frac{4\cdot 3!}{2^7} + \frac{7(2^4+2)}{2^{11}} \,<\, \frac{1}{4},
\end{align*}
since $\eta < 2^{-11}$, as required.
\end{proof}


Finally, in order to show that $\bigcup_{I \in \cI}I$ is disjoint from the set
$(\pi/2)\Z + [-100\pi/n,100\pi/n]$, we will need the following simple lemma.


\begin{lemma}\label{lem:origin}
 If\/ $|x| \le 1/8$ or\/ $|x - T\pi| \le 1/8$, then $\Re\big(H(x)\big)\ge 1/2$.
\end{lemma}
\begin{proof}
We will use the following facts (cf.~\cite[Theorem~5]{BC}), which can be easily verified by induction:
for every $t \ge 0$,
\[
 P_{2t}(1) = P_{2t}(-1) = Q_{2t}(1) = -Q_{2t}(-1) = 2^t,
\]
and
\[
 P_{2t+1}(1) = Q_{2t+1}(-1) = 2^{t+1},\qquad P_{2t+1}(-1) = Q_{2t+1}(1) = 0.
\]
Since $t$ is odd, it follows that
$$\Re\big(H(0)\big) = 2^{-(t+1)/2}\big(P_t(1) +Q_t(1) \big)=1$$
and
$$\Re\big(H(T\pi)\big) = 2^{-(t+1)/2}\big(P_t(-1)+Q_t(-1)\big)=1.$$
Now, by Lemma~\ref{lem:H_deriv_ubounds} we have $|H'(x)|\le 4$
for every $x \in \R$, and so
\[
 \Re\big(H(x)\big) \ge 1-4|x| \ge \frac{1}{2}
\]
for all $x \in \R$ with $|x|\le 1/8$. A similar argument works for those $x$ near $T \pi$.
\end{proof}

\begin{remark}
Note that $x = T\pi$ corresponds to $\theta = \pi/2$ in the cosine polynomial $c(\theta)$.
The reader may have noticed that we do not necessarily need the cosine polynomial
to be large at this point, as the sine polynomial can be large there. However, for technical
reasons, this will be useful later on, in the proof of Lemma~\ref{lem:bounding:s:hat}.
\end{remark}

We are finally ready to prove Theorem~\ref{thm:cosine}.

\begin{proof}[Proof of Theorem~\ref{thm:cosine}]
Let $c(\theta)$ be the cosine polynomial defined in~\eqref{def:c}, and recall that $\supp(c) = C$,
that $\eps_k \in \{-1,1\}$ for every $k \in C$, and that $|c(\theta)|\le \sqrt{n}$ for every
$\theta \in \R/2\pi\Z$, by Lemma~\ref{lem:even_ubounds}. We will show that there exists a suitable
and well-separated collection $\cI$ of disjoint intervals in $\R/2\pi\Z$ such that
$|c(\theta)| \ge \delta \sqrt{n}$ for all $\theta \notin \bigcup_{I \in \cI}I$.

To prove this, set $\eta := 2T\pi / n$, and note that $\eta < \pi\gamma < 2^{-11}$.
Partition $\R/4T\pi\Z=\R/2n\eta\Z$ into $2n$ intervals $I_j := [j\eta,(j+1)\eta]$, each of length~$\eta$,
and say that an interval $I_j$ is {\em good\/} if
\[
 \big| \Re\big( H(x) \big) \big| \ge \frac{\eta^3}{2^7}
\]
for all $x \in I_j$. Let $\cJ'$ be the collection of maximal unions of consecutive good
intervals~$I_j$, and let $\cI'$ be the collection of remaining intervals (i.e., maximal
unions of consecutive bad intervals). Thus $\cI'$ and $\cJ'$ form interleaving collections
of intervals decomposing $\R/4T\pi\Z$. Scaling from $x$ to $\theta = x/2T$ 
gives corresponding collections of intervals $\cI$ and $\cJ$; we claim that $\cI$ is the required
suitable and well-separated collection.

First, to see that $\cI$ is suitable, note that each interval $I_j$ (and hence each $I \in \cI'$)
starts and ends at a multiple of $\eta=2T\pi/n$. Hence after scaling, each $I\in\cI$ starts and ends
at points of $\frac{\pi}{n}\Z$. The set $\cI$ is invariant under the maps
$\theta \mapsto \pi \pm \theta$ by the symmetries of the function $\cos(k\theta)$ when $k\in C\subseteq2\Z$.
To see that $|\cI| \le 4\gamma n$, note that since a cosine polynomial
of degree $d$ has at most $2d$ roots in its period, there are at most $4(2T + 2^t - 1) = 4\gamma n$
values of $x \in \R/4T\pi\Z$ where $\Re(H(x)) = 2^{-7} \eta^3$, and the same bound on the number where
$\Re(H(x)) = -2^{-7} \eta^3$. Since each $I \in \cI'$ must contain at least two such points
(counted with multiplicity), we have $|\cI| = |\cI'| \le 4\gamma n$, as required.

Next, let us show that $\cI$ is well-separated. Recall first that, by Lemma~\ref{lem:intervals:nice},
any set of 7 consecutive intervals $I_j$ must contain a good interval.
Thus $|I|\le 6\eta$ for each $I\in\cI'$, and so $|I|\le 6\pi/n$ for each $I\in\cI$.
Now, $d(I,J) \ge \pi/n$ for distinct $I,J \in \cI$ by construction, and the sets $[-100\eta,100\eta]$
and $T\pi + [-100\eta,100\eta]$ are each contained in an element of $\cJ'$ by Lemma~\ref{lem:origin},
since $2^{-7}\eta^3<1/2$ and $100\eta<1/8$. Scaling down, it follows that $\bigcup_{I \in \cI}I$
is disjoint from the set $(\pi/2)\Z + [-100\pi/n,100\pi/n]$, as required.

Finally, recalling that $\eta = 2T\pi / n$, $\gamma n=2T+2^t-1$, $T=2^{t+10}$,
and that $|\Re(H(x))|\ge 2^{-7} \eta^3$ for each $x \in J \in \cJ'$, it follows that
\[
 |c(\theta)| \ge 2^{(t+1)/2}\cdot 2^{-7}\eta^3=2^{-12}\pi^3(2T)^{7/2}/n^3
 \ge 2^{-8}\gamma^{7/2}\sqrt{n}=\delta \sqrt{n}
\]
for every $\theta \notin \bigcup_{I \in \cI}I$, as required.
\end{proof}

\section{Minimising Discrepancy}\label{brownian:sec}

In this section we recall the main `partial colouring' lemma of Spencer~\cite{S85} (whose proof, as noted in the introduction, was based on a technique of Beck~\cite{B81}), which will play an important role in the proof of Theorem~\ref{thm:sine}. In particular, we will use the results of this section both to choose in which direction we should `push' the sine polynomial on each interval $I \in \cI$, and to show that we can choose $\eps_k \in \{-1,1\}$ so that it is pushed (roughly) the correct distance. The following convenient variant of Spencer's theorem was proved by Lovett and Meka~\cite[Theorem~4]{LM}\footnote{The theorem as stated in~\cite{LM} only insists that $|x_i|\ge 1-\delta$ for at least $n/2$ indices, due to the requirement that a fast algorithm exists. However, it is clear by continuity that we can take $\delta=0$ if we are only interested in an `existence proof'.}, who also gave a beautiful polynomial-time randomised algorithm for finding a colouring with small discrepancy. 

\begin{theorem}[Main Partial Colouring Lemma]\label{thm:partial_colouring}
 Let\/ $v_1,\dots,v_m \in \R^n$ and\/ $x_0\in[-1,1]^n$. If\/ $c_1,\dots,c_m \ge 0$ are such that
 \[
  \sum_{j=1}^m\exp\big(-c_j^2/16\big) \le \frac{n}{16},
 \]
 then there exists an\/ $x \in [-1,1]^n$ such that
 \[
  |\langle x-x_0,v_j\rangle|\le c_j\|v_j\|_2
 \]
 for every\/ $j \in [m]$, and moreover\/ $x_i \in \{-1,1\}$
 for at least\/ $n/2$ indices\/ $i \in [n]$.
\end{theorem}

We will in fact use the following corollary of Theorem~\ref{thm:partial_colouring}.

\begin{corollary}\label{cor:colouring}
 Let\/ $v_1,\dots,v_m \in \R^n$ and\/ $x_0\in[-1,1]^n$. If\/ $c_1,\dots,c_m \ge 0$ are such that
 \begin{equation}\label{e:c_j_bound}
  \sum_{j=1}^m\exp\big(-c_j^2/14^2\big) \le \frac{n}{16},
 \end{equation}
 then there exists an\/ $x\in\{-1,1\}^n$ such that
 \[
  |\langle x-x_0,v_j\rangle|\le (c_j+30)\sqrt{n} \cdot \|v_j\|_\infty
 \]
 for every\/ $j \in [m]$.
\end{corollary}

\begin{proof}
We prove Corollary~\ref{cor:colouring} by induction on~$n$. Note first that the result is trivial for all $n \le 900$,
since we can choose $x \in \{-1,1\}^n$ with $\| x - x_0 \|_\infty \le 1$, and for such a vector we have $|\langle x - x_0, v_j \rangle| \le n \cdot \|v_j\|_\infty \le 30\sqrt{n} \cdot \|v_j\|_\infty$. 

For $n > 900$, we apply Theorem~\ref{thm:partial_colouring} with constants $b_j := 2c_j/7$,  noting that
\[
 \sum_{j=1}^m \exp\big(-b_j^2/16\big)=\sum_{j=1}^m\exp\big(-c_j^2/14^2\big)\le \frac{n}{16}.
\]
We obtain a vector $y \in [-1,1]^n$, with
\[
 |\langle y - x_0, v_j \rangle| \le b_j \|v_j\|_2 \le b_j\sqrt{n} \cdot \|v_j\|_\infty
\]
for every $j \in [m]$, such that $y_i \in \{-1,1\}$ for at least $n/2$ indices $i \in [n]$.

Now, let $U \subseteq [n]$ be a set of size $\lceil n/2 \rceil$ such that $y_i \in \{-1,1\}$
for every $i \in U$, and set $W := [n] \setminus U$. For each $j \in [m]$, define a constant $a_j \ge 0$ so that
\[
 a_j^2 := c_j^2 + 14^2 \log\bigg( \frac{n}{\lfloor n/2 \rfloor} \bigg),
\]
and observe that
\[
 \sum_{j = 1}^m \exp\big(-a_j^2/14^2\big) \le \frac{\lfloor n/2\rfloor}{16} = \frac{|W|}{16},
\]
and that $a_j \le c_j + 12$, since $14^2\log(n/\lfloor n/2\rfloor)<196\log 2.01<12^2$ for $n>900$.

Let $\pi \colon \R^n\to\R^W$ be projection onto the coordinates of~$W$.
By the induction hypothesis, we obtain a vector $z \in \{-1,1\}^W$ with
\[
 |\langle z - \pi(y),\pi(v_j) \rangle| \le (a_j + 30)\sqrt{|W|} \cdot \|\pi(v_j)\|_\infty
 \le (a_j+30)\sqrt{n/2} \cdot \|v_j\|_\infty.
\]
Now, define $x \in \{-1,1\}^n$ by setting $x_i := y_i$ for $i \in U$ and $\pi(x) = z$,
and observe that
\begin{align*}
 |\langle x-x_0,v_j\rangle|
 & \le |\langle y - x_0,v_j \rangle| + |\langle z - \pi(y), \pi(v_j)\rangle|\\
 & \le \big( b_j + (a_j + 30)/\sqrt{2} \big) \sqrt{n} \cdot \|v_j\|_\infty\\
 & \le \bigg( \frac{2c_j}{7} + \frac{c_j + 42}{\sqrt{2}} \bigg) \sqrt{n} \cdot \|v_j\|_\infty\\
 & \le (c_j+30)\sqrt{n}\cdot\|v_j\|_\infty,
\end{align*}
as required, since $b_j = 2c_j/7$ and $a_j \le c_j + 12$. This completes the induction step.
\end{proof}

\begin{remark}
 The result is stated in terms of the $\ell^\infty$-norms $\|v_j\|_\infty$ because we cannot control
 the decrease in $\|v_j\|_2$ when we discard half of the coordinates.
\end{remark}

\begin{remark}
 It is important for our application that $m$ can be much larger than~$n$,
 and that the only restriction on $m$ occurs via the condition~\eqref{e:c_j_bound}.
 In particular, we will later apply Corollary~\ref{cor:colouring} with $m$ very large,
 but with the $c_j$ increasing sufficiently rapidly so that~\eqref{e:c_j_bound} still holds.
\end{remark}

\section{The odd sine polynomial}\label{sine:sec}

The aim of this section is to prove Theorem~\ref{thm:sine}. Let $\cI$ be a collection of suitable
well-separated intervals, and recall from Definition~\ref{def:wellseparated} that $|\cI| = 4N$ for some $N \le \gamma n$, and that $\cI$ is invariant under the maps $\theta \mapsto \pi \pm \theta$. The collection $\cI$ is therefore uniquely determined by the set $\cI_0 \subseteq \cI$ of $N$ intervals that lie in $[0,\pi/2]$ (since no $I \in \cI$ contains 0 or $\pi/2$).

As described in Section~\ref{strategy:sec}, our aim is to `push' the sine polynomial away from zero (in either the positive or negative direction) on each interval in $\cI$. Let us say that a colouring $\alpha \colon \cI \to \{-1,1\}$ is \emph{symmetric\/} if $\alpha(I') = \alpha(I)$ whenever $I' = \pi - I$, and $\alpha(I') = -\alpha(I)$ whenever $I' = \pi + I$. Note that if $\alpha$ is symmetric, then it is uniquely determined by its values on the set $\cI_0$. Finally, recall that $S_o=\{1,3,5,\dots,2n-1\}$, and set $K := 2^7$. 


\begin{definition}
 Given a colouring $\alpha \colon \cI \to \{-1,1\}$, we define
 $g_{\alpha} \colon \R/2\pi\Z \to \{-1,0,1\}$ by
 \[
  g_\alpha(\theta):=\sum_{I \in \cI} \alpha(I) \id[\theta \in I].
 \]
 We also define a vector $\hat\eps = (\hat\eps_1,\hat\eps_3,\ldots,\hat\eps_{2n-1}) \in \R^{S_o}$ by setting  
 \[
  \hat\eps_j := K \sqrt{n} \int_{-\pi}^{\pi} g_{\alpha}(\theta)\sin(j\theta)\,d\theta,
 \]
 for each $j \in S_o$. 
\end{definition}

\begin{remark}
 By Fourier inversion, one would expect the function
 $\hat s_\alpha(\theta):=\sum_{j \in S_o}\hat\eps_j \sin(j\theta)$
 to approximate $\pi K\sqrt{n}\,g_{\alpha}(\theta)$; in particular, it should be large on the intervals $I \in \cI$. 
 We will prove in Lemma~\ref{lem:bounding:s:hat}, below, that this is indeed the case. 
\end{remark}

We will use $\hat\eps$ as the starting point of an application of Corollary~\ref{cor:colouring},
so we need $|\hat\eps_j| \le 1$ for all $j \in S_o$. The following lemma, which we also prove
using Corollary~\ref{cor:colouring}, shows that, since we chose $\gamma$ sufficiently small,
we can choose the colouring $\alpha$ so that this is the case.

\begin{lemma}\label{lem:choosing:alpha}
 There exists a symmetric colouring\/ $\alpha \colon \cI \to \{-1,1\}$ such that\/ $\hat\eps \in [-1,1]^{S_o}$.
\end{lemma}
\begin{proof}
Write $\cI_0=\{I_1,\dots,I_N\}$ and recall that this collection determines~$\cI$.
Now, for each $j \in [n]$, define a vector $v_j \in \R^N$ by setting
\[
 (v_j)_i := 4K\sqrt{n} \int_{I_i} \sin\big((2j-1)\theta\big)\,d\theta
\]
for each $i \in [N]$, and observe that, for each $j \in [n]$, we have
\[
 \hat\eps_{2j-1} = K \sqrt{n} \int_{-\pi}^{\pi} g_{\alpha}(\theta)\sin\big( (2j-1)\theta \big)\,d\theta
 =\sum_{i=1}^N \alpha(I_i) (v_j)_i,
\]
by the symmetry conditions on both $\alpha$ and~$\cI$. Our task is therefore to find a vector
$x \in \{-1,1\}^N$ such that $|\langle x,v_j \rangle|\le 1$ for all $j \in [n]$. Indeed, we will then
be able to set $\alpha(I_i)=x_i$ for each $i \in [N]$, and deduce that $|\hat\eps_k|\le1$ for all $k\in S_o$.

We do so by applying Corollary~\ref{cor:colouring} with $x_0 := 0$ and
$c_j := 14 \sqrt{\log(16n / N)}$ for each $j \in [n]$. Noting that~\eqref{e:c_j_bound} is satisfied,
it follows from Corollary~\ref{cor:colouring} that there exists an $x\in\{-1,1\}^N$ such that
\[
 |\langle x,v_j\rangle| \le \big(c_j+30\big) \sqrt{N} \cdot \|v_j\|_\infty.
\]
Now, since $\cI$ is well-separated, by Definition~\ref{def:wellseparated}$(d)$ we have
\[
 |(v_j)_i| \le 4K\sqrt{n} \cdot |I_i| \le \frac{24\pi K}{\sqrt{n}}
\]
for every $i \in [N]$ and $j \in [n]$. It follows that
\[
 |\langle x,v_j\rangle| \le \big(14\sqrt{\log(16n/N)}+30\big)\sqrt{N/n}\cdot 24\pi K.
\]
Note that the right hand side is an increasing function of $N$ for $N/n\le \gamma < 1$ and so
\[
 |\langle x,v_j\rangle| \le \big(14\sqrt{\log(16/\gamma)}+30\big)\sqrt{\gamma}\cdot 24\pi K\le 1,
\]
where the last inequality follows from our choice of $K=2^7$ and the inequality $\gamma\le 2^{-40}$.
\end{proof}

For the rest of the proof fix this colouring $\alpha$ (and hence also the vector~$\hat\eps$).
Recall that our aim is to choose a colouring $\eps \colon S_o \to \{-1,1\}$ so that the conclusion
of Theorem~\ref{thm:sine} holds. Given such a colouring, define
\[
 s_o(\theta) := \sum_{j \in S_o} \eps_j \sin(j\theta) \qquad\qquad \text{and} \qquad\qquad
 \hat s_\alpha(\theta) := \sum_{j \in S_o} \hat\eps_j \sin(j\theta).
\]
Our aim is to choose the $\eps_j$ so that $|s_o(\theta)-\hat s_\alpha(\theta)|$ is uniformly bounded
for all $\theta\in\R/2\pi\Z$ (see Lemma~\ref{lem:choosing:epsilon}, below). A na{\"\i}ve approach to controlling this difference on a sufficiently
dense set of points would require imposing more constraints (with smaller values of $c_j$) than can be handled by Corollary~\ref{cor:colouring}. Instead we shall
place constraints on the differences $|s_o^{(\ell)}(\theta) - \hat s_\alpha^{(\ell)}(\theta)|$ of the $\ell$th
derivatives for each $\ell \ge 0$, but at many fewer values of~$\theta$, and then use Taylor's Theorem
to bound $|s_o(\theta)-\hat s_\alpha(\theta)|$ at all other points. The advantage of this approach
is that the constraints we need on the higher derivatives become rapidly weaker as $\ell$
increases, and in particular can be chosen so that~\eqref{e:c_j_bound} is satisfied.

Note that it is enough to bound $|s_o(\theta)-\hat s_\alpha(\theta)|$ on $[0,\frac{\pi}{2}]$ as
both $s_o(\theta)$ and $\hat s_\alpha(\theta)$ have the same symmetries under $\theta\mapsto \pi\pm\theta$.
Set $M:=16n$ and let $\theta_k:=\frac{(2k-1)\pi}{4M}$ for $k=1,\dots,M$.
Then for any point $\theta \in [0,\frac{\pi}{2}]$
there exists $k\in[M]$ such that $|\theta-\theta_k|\le \frac{\pi}{4M}= 2^{-6} \pi / n$.
By Taylor's Theorem (and the fact that all sine polynomials are entire functions so their Taylor expansions converge), we have
\begin{equation}\label{eq:Taylor}
 s_o(\theta)-\hat s_\alpha(\theta) = \sum_{\ell=0}^\infty
 \big( s_o^{(\ell)}(\theta_k)-\hat s_\alpha^{(\ell)}(\theta_k)\big)\frac{(\theta-\theta_k)^\ell}{\ell!}.
\end{equation}
We will bound the absolute value of the right-hand side using Corollary~\ref{cor:colouring}.

\pagebreak

\begin{lemma}\label{lem:sin_derivatives}
 There exists a colouring\/ $\eps \colon S_o \to \{-1,1\}$ such that
 \[
  \big|s_o^{(\ell)}(\theta_k)-\hat s_\alpha^{(\ell)}(\theta_k)\big| \le (65+2\ell)\sqrt{n} \cdot (2n)^\ell
 \]
 for every\/ $k \in [M]$ and\/ $\ell \ge 0$.
\end{lemma}

\begin{proof}
For each $k \in [M]$ and $\ell \ge 0$, define a vector $v_{(k,\ell)}\in \R^n$ by setting
\[
 (v_{(k,\ell)})_j = \frac{d^\ell}{d\theta^\ell} \sin\big((2j-1)\theta\big) \big|_{\theta = \theta_k}
\]
for each $j \in [n]$, and observe that
\[
 s_o^{(\ell)}(\theta_k) - \hat s_\alpha^{(\ell)}(\theta_k)
 = \sum_{j=1}^n \big(\eps_{2j-1} - \hat\eps_{2j-1}\big) (v_{(k,\ell)})_j
 = \langle \eps-\hat\eps, v_{(k,\ell)}\rangle,
\]
where we consider $\eps - \hat\eps$ and $v_{(k,\ell)}$ as vectors in $\R^{S_o}$.

We apply Corollary~\ref{cor:colouring} with $x_0 := \hat\eps$ and
$c_{(k,\ell)}=14\sqrt{(9+\ell)\log 2}$. Observe that
\[
 \sum_{k=1}^M \sum_{\ell=0}^\infty \exp\big(-c_{(k,\ell)}^2/14^2\big)
 =\sum_{k=1}^M \sum_{\ell=0}^\infty 2^{-(9+\ell)} = M \cdot 2^{-8} = \frac{n}{16},
\]
and so~\eqref{e:c_j_bound} is satisfied. It follows\footnote{Note that we
appear to be applying Corollary~\ref{cor:colouring} with an infinite number of constraints,
but in fact only finitely many of them are needed as the constraints vacuously hold when $\ell \ge n$.}
from Corollary~\ref{cor:colouring} that there exists an $\eps \in \{-1,1\}^n$ such that
\[
 |\langle \eps-\hat\eps, v_{(k,\ell)}\rangle|
 \le \big( c_{(k,\ell)}+30 \big) \sqrt{n} \cdot \|v_{(k,\ell)}\|_\infty
\]
for every $k \in [M]$ and $\ell \ge 0$. Now, observe that
\[
 \|v_{(k,\ell)}\|_\infty \le (2n)^\ell,
\]
and that $14^2(9+\ell)\log 2\le 35^2+140\ell \le (35+2\ell)^2$, so
\[
 c_{(k,\ell)}+30 \le 65+2\ell.
\]
Combining these bounds, we obtain
\[
 \big|s_o^{(\ell)}(\theta_k)-\hat s_\alpha^{(\ell)}(\theta_k)\big|
 = |\langle \eps-\hat\eps, v_{(k,\ell)}\rangle| \le (65+2\ell) \sqrt{n} \cdot (2n)^\ell
\]
for every $k \in [M]$ and $\ell \ge 0$, as required.
\end{proof}

The following bound on the magnitude of $s_o(\theta) - \hat s_\alpha(\theta)$ is a straightforward consequence.

\pagebreak

\begin{lemma}\label{lem:choosing:epsilon}
 There exists a colouring\/ $\eps \colon S_o \to \{-1,1\}$ such that
 \[
  |s_o(\theta)-\hat s_\alpha(\theta)| \le 72\sqrt{n}
 \]
 for every\/ $\theta \in \R$.
\end{lemma}

\begin{proof}
Let us assume (without loss of generality) that $\theta \in [0,\frac{\pi}{2}]$, and let $k\in[M]$ be such that $|\theta-\theta_k| \le 2^{-6} \pi / n$. By~\eqref{eq:Taylor} and Lemma~\ref{lem:sin_derivatives}, we have
\[
 \big|s_o(\theta)-\hat s_\alpha(\theta)\big|
 \le \sum_{\ell=0}^\infty \big| s_o^{(\ell)}(\theta_k) - \hat s_\alpha^{(\ell)}(\theta_k) \big|
 \frac{(2^{-6} \pi / n)^\ell}{\ell!}
 \le \sum_{\ell = 0}^\infty (65+2\ell)\sqrt{n} \cdot \frac{(2^{-5} \pi)^\ell}{\ell!}.
\]
Now simply observe that
\[
 \sum_{\ell = 0}^\infty (65+2\ell) \frac{(2^{-5}\pi)^\ell}{\ell!} = (65 + 2^{-4} \pi) e^{2^{-5} \pi} \le 72,
\]
and the lemma follows.
\end{proof}

We will prove that the conclusion of Theorem~\ref{thm:sine} holds for the colouring $\eps$ given by
Lemma~\ref{lem:choosing:epsilon}. To deduce this, it will suffice to show that $\hat s_\alpha(\theta)$ approximates the
step function $\pi K\sqrt{n} \cdot g_\alpha(\theta)$ sufficiently well, and in particular that it is large
on each interval~$I \in \cI$.

\begin{lemma}\label{lem:bounding:s:hat}
 For every\/ $\theta \in \bigcup_{I \in \cI}I$, we have
 \[
  |\hat s_\alpha(\theta)| \ge \frac{2K\sqrt{n}}{3}.
 \]
 Moreover, $|\hat s_\alpha(\theta)|\le 5K\sqrt{n}$ for every\/ $\theta \in \R$.
\end{lemma}

The proof of Lemma~\ref{lem:bounding:s:hat} follows from a standard (but somewhat technical) calculation, and to simplify things slightly we will find it convenient to renormalise, by defining
\[
 \tg_\alpha(\theta) := (K\sqrt{n})^{-1} \hat s_\alpha(\theta).
\]
Fix $\theta_0 \in \R$, and observe that, by the symmetry conditions on both $\alpha$ and~$\cI$, we have
\begin{align}
 \tg_\alpha(\theta_0)
 & \,=\, \sum_{j=0}^{n-1} \sin\big( (2j+1)\theta_0 \big)
   \int_{-\pi}^{\pi} g_\alpha(\theta) \sin\big( (2j+1)\theta \big)\,d\theta \nonumber\\
 & \,=\, 4\int_0^{\pi/2}g_\alpha(\theta)\sum_{j=0}^{n-1}
   \sin\big( (2j+1)\theta_0 \big) \sin\big( (2j+1)\theta \big)\,d\theta.\label{eq:stilde:firststep}
\end{align}
We can now use the following standard trigonometric fact.

\begin{observation}\label{obs:trig}
\[
 4 \sum_{j=0}^{n-1} \sin\big( (2j+1)\theta_0 \big) \sin\big( (2j+1)\theta \big)
 = \frac{\sin\big( 2n(\theta-\theta_0) \big)}{\sin(\theta-\theta_0)}
 - \frac{\sin\big( 2n(\theta+\theta_0) \big)}{\sin(\theta+\theta_0)}.
\]
\end{observation}

\begin{proof}
Simply note that both sides are equal to
\[
 2 \sum_{j=0}^{n-1} \Big( \cos\big( (2j+1)(\theta-\theta_0) \big) - \cos\big( (2j+1)(\theta+\theta_0) \big) \Big),
\]
using the addition formulae for $\sin(\alpha\pm\beta)$ and $\cos(\alpha\pm\beta)$
and the telescoping series
\begin{align*}
 \sin(2n\varphi)&=\sum_{j=0}^{n-1} \Big( \sin\big( (2j+1)\varphi + \varphi \big) - \sin\big( (2j+1)\varphi - \varphi \big) \Big)\\
  &=\sum_{j=0}^{n-1} 2 \cos\big( (2j+1)\varphi \big) \sin(\varphi)
\end{align*}
for $\varphi=\theta\pm\theta_0$.
\end{proof}

Combining~\eqref{eq:stilde:firststep} and Observation~\ref{obs:trig}, and recalling the definition of $g_{\alpha}(\theta)$, it follows that
\begin{equation}\label{eq:stilde:twostep}
 \tg_\alpha(\theta_0) = \sum_{I \in \cI_0} \alpha(I) \int_I \bigg( \frac{\sin\big( 2n(\theta-\theta_0) \big)}{\sin(\theta-\theta_0)}
 - \frac{\sin\big( 2n(\theta+\theta_0) \big)}{\sin(\theta+\theta_0)} \bigg)\,d\theta.
\end{equation}
Before bounding the right-hand side of~\eqref{eq:stilde:twostep}, let us briefly discuss what is going on. Let $\theta_0 \in [0,\pi/2]$, and recall from Definition~\ref{def:wellseparated}(f) that no $I \in \cI_0$ contains any point close to $0$ or $\pi/2$. It follows that the integrand in~\eqref{eq:stilde:twostep} behaves roughly like a point mass placed at $\theta = \theta_0$, and hence $\tg_\alpha(\theta_0)$ should be approximately $\alpha(I)$ when $\theta_0 \in I$, and small otherwise. 

To make this rigorous, we will show that the integral of the first term over the interval $I \in \cI_0$ containing $\theta_0$ (if such an interval exists) is of order 1, and that the integral over the remaining intervals (and over the second term) is smaller. This will follow via a straightforward calculation from the fact that the endpoints of each interval in $\cI$ lie in $\frac{\pi}{n} \Z$. 



Instead of approximating the integral for intervals close to $\theta_0$ directly, we will instead compare it to a slightly simpler `sine integral', which we bound in the following lemma.

\begin{lemma}\label{lem:si}
 Let\/ $I \in \cI$ and let\/ $\theta_0 \in \R$.
 \begin{itemize}
  \item[$(a)$] If\/ $\theta_0\in I$, then
   \[
     \frac{4}{3} \,\le\, \int_I\frac{\sin\big( 2n(\theta-\theta_0) \big)}{\theta-\theta_0}\,d\theta \,\le\, 4.
   \]
  \item[$(b)$]  If\/ $\theta_0\notin I$ then
   \[
    \!\!-1 \,\le\, \int_I\frac{\sin\big( 2n(\theta-\theta_0) \big)}{\theta-\theta_0}\,d\theta \,\le\, 2.
   \]
 \end{itemize}
\end{lemma}
\begin{proof}
Recall from Definition~\ref{def:wellseparated} that the endpoints of $I$ are in\/ $\frac{\pi}{n}\Z$,
and let $I = [a\pi/n,b\pi/n]$, where $a,b \in \Z$ with $a < b$. Substituting $x=2n(\theta-\theta_0)$
gives us the integral
\[
 f(\theta_0):=\int_I\frac{\sin\big( 2n(\theta-\theta_0) \big)}{\theta-\theta_0}\,d\theta
 = \int_{2a\pi-2n\theta_0}^{2b\pi-2n\theta_0}\frac{\sin x}{x}\,dx,
\]
and we note that
\[
 f'(\theta_0) = (-2n)\bigg( \frac{\sin(2b\pi-2n\theta_0)}{2b\pi-2n\theta_0}
  -\frac{\sin(2a\pi-2n\theta_0)}{2a\pi-2n\theta_0} \bigg)
  = \frac{4\pi n(a-b)\sin(2n\theta_0)}{(2a\pi-2n\theta_0)(2b\pi-2n\theta_0)},
\]
since $a,b \in \Z$, so $\sin(2a\pi-2n\theta_0) = \sin(2b\pi-2n\theta_0) = - \sin(2n\theta_0)$.
Since $a \ne b$, it follows that the extremal values of $f(\theta_0)$ can occur only
when $\sin(2n\theta_0) = 0$, i.e., when $2n\theta_0\in\pi\Z$. These extremal values
must therefore be of the form
\[
 u(\ell) + u(\ell+1) + \dots + u(\ell + 2(b-a) - 1)
\]
for some $\ell \in \Z$, where
\[
 u(j) := \int_{j\pi}^{(j+1)\pi}\frac{\sin x}{x}\,d\theta.
\]
We claim first that if $\theta_0\in I$, then
\[
 \int_{0}^{2\pi}\frac{\sin x}{x}\,d\theta \le f(\theta_0) \le \int_{-\pi}^{\pi}\frac{\sin x}{x}\,d\theta.
\]
Indeed, if $\theta_0\in I$ then $2a\pi \le 2\theta_0 n \le 2b\pi$, and so
$\ell \le 0 \le \ell + 2(b-a)$. Note also that
\[
 u(2j) > 0, \qquad u(2j+1) < 0 \qquad \text{and} \qquad u(-j) = u(j-1)
\]
for every non-negative $j \in \Z$, and moreover
\[
 u(2j-1) + u(2j) < 0 < u(2j) + u(2j+1)
\]
for every $j \ge1$. It follows that the maxima of $f(\theta_0)$ are at most
$u(-1) + u(0)$, and the minima are at least $u(0) + u(1)$, as claimed.
Similarly, if $\theta_0 \notin I$ then without loss of generality we have $\ell \ge 0$,
and by the same argument as above we have
\[
 \int_\pi^{2\pi} \frac{\sin x}{x}\,d\theta \le f(\theta_0) \le \int_0^\pi \frac{\sin x}{x}\,d\theta.
\]
It is now straightforward to obtain the claimed bounds by numerical integration.
\end{proof}

We will use the following simple lemma to bound the integrals
in the proof of Lemma~\ref{lem:bounding:s:hat}.

\begin{lemma}\label{lem:periodic}
 If\/ $h\colon[a,b]\to\R$ is a monotonic function and\/ $b-a\in\tfrac{\pi}{n}\Z$, then
 \[
  \bigg| \int_a^b h(\theta)\sin(2n\theta)\,d\theta \bigg| \le \frac{|h(b) - h(a)|}{n}.
 \]
\end{lemma}
\begin{proof}
Assume without loss of generality that $h$ is increasing, and suppose first that $b=a+\frac{\pi}{n}$.
Since $\sin(x+\pi)=-\sin(x)$ we have
\[
 \int_a^{a+\frac{\pi}{n}} h(\theta)\sin(2n\theta)\,d\theta
 =\int_a^{a+\frac{\pi}{2n}} \big( h(\theta) - h(\theta+\tfrac{\pi}{2n}) \big) \sin(2n\theta)\,d\theta,
\]
and therefore, since $h$ is increasing,
\[
 \bigg| \int_a^{a+\frac{\pi}{n}} h(\theta)\sin(2n\theta)\,d\theta \bigg|
 \,\le\, \big( h(b) - h(a) \big) \int_a^{a + \frac{\pi}{2n}}|\sin(2n\theta)|\,d\theta
 \,=\, \frac{h(b) - h(a)}{n},
\]
as required. To deduce the general case, simply split the interval $[a,b]$ into sub-intervals
of length $\frac{\pi}{n}$ and use the triangle inequality.
\end{proof}

We are now ready to prove Lemma~\ref{lem:bounding:s:hat}.

\begin{proof}[Proof of Lemma~\ref{lem:bounding:s:hat}]
Recall that it is enough to prove the bounds when $\theta=\theta_0\in[0,\pi/2]$, and that $\cI_0=\{I\in\cI:I\subseteq[0,\pi/2]\}$. 
By~\eqref{eq:stilde:twostep}, we have
\begin{equation}\label{eq:stilde:twostep:again}
 \tg_\alpha(\theta_0) = \sum_{I \in \cI_0} \alpha(I) \int_I \bigg( \frac{\sin\big( 2n(\theta-\theta_0) \big)}{\sin(\theta-\theta_0)}
 - \frac{\sin\big( 2n(\theta+\theta_0) \big)}{\sin(\theta+\theta_0)} \bigg)\,d\theta
\end{equation}
for every $\theta_0 \in [0,\pi/2]$. We will deal with the second term first.

\medskip
\noindent \textbf{Claim 1:}
$\ds\sum_{I \in \cI_0} \bigg| \int_I \frac{\sin(2n(\theta+\theta_0))}{\sin(\theta+\theta_0)}\,d\theta \bigg|
\,\le\, \ds\frac{1}{50\pi} + \frac{O(1)}{n}$.

\begin{proof}[Proof of Claim~1]
Let $I \in \cI_0$, and suppose first that $\sin \theta$ is monotonic on $I + \theta_0$. By Lemma~\ref{lem:periodic},
applied with $h(\theta)=1/\sin\theta$, we have
\[
 \bigg| \int_I \frac{\sin(2n(\theta+\theta_0))}{\sin(\theta+\theta_0)}\,d\theta \bigg|
 \,\le\, \frac{1}{n} \bigg(\max_{\theta\in I} \frac{1}{\sin(\theta+\theta_0)}
  \,-\, \min_{\theta\in I} \frac{1}{\sin(\theta+\theta_0)} \bigg),
\]
since, by Definition~\ref{def:wellseparated}, the endpoints of $I$ are in\/ $\frac{\pi}{n}\Z$.
If $\sin \theta$ is not monotonic on $I + \theta_0$, then we instead use the trivial bound
\[
 \bigg| \int_I \frac{\sin(2n(\theta+\theta_0))}{\sin(\theta+\theta_0)}\,d\theta \bigg|
 \,\le\, |I| \cdot \max_{\theta\in I} \frac{1}{\sin(\theta+\theta_0)}
 \,=\, \frac{O(1)}{n},
\]
where the final inequality holds since $|I| = O(1/n)$, by Definition~\ref{def:wellseparated},
and hence (since $\sin \theta$ is not monotonic on $I + \theta_0\subseteq[0,\pi]$)
we have $\sin(\theta+\theta_0) > 1/2$ for all $\theta \in I$.

Now, summing over intervals $I \in \cI_0$, and partitioning into three classes according
to whether $\sin \theta$ is increasing, decreasing, or neither on $I + \theta_0$, we obtain
two alternating sums that are both bounded by their maximum terms, and possibly one additional term
(for which we use the trivial bound). Recalling from Definition~\ref{def:wellseparated}
that $\bigcup_{I \in \cI}I$ is disjoint from the set $(\pi/2)\Z + [-100\pi/n,100\pi/n]$, we obtain
\[
\ds\sum_{I \in \cI_0} \bigg| \int_I \frac{\sin(2n(\theta+\theta_0))}{\sin(\theta+\theta_0)}\,d\theta \bigg|
 \,\le\, \frac{2}{n\sin(100\pi/n)} + \frac{O(1)}{n}
 \,=\, \frac{1}{50\pi} + \frac{O(1)}{n},
\]
as claimed.
\end{proof}

The next claim will allow us to replace the first term in~\eqref{eq:stilde:twostep:again}
by the integral in Lemma~\ref{lem:si}.

\medskip
\noindent \textbf{Claim 2:}
$\ds\sum_{I \in \cI_0} \bigg| \int_I \frac{\sin\big( 2n(\theta-\theta_0) \big)}{\sin(\theta-\theta_0)} \,
- \, \frac{\sin\big( 2n(\theta-\theta_0) \big)}{\theta - \theta_0} \,d\theta \bigg|
= \ds\frac{O(1)}{n}$.

\begin{proof}[Proof of Claim~2]
We again apply Lemma~\ref{lem:periodic}, this time with $h(\theta) = \frac{1}{\sin\theta}-\frac{1}{\theta}$, which
is increasing on $[-\pi/2,\pi/2]$, to give
\[
 \bigg| \int_I\frac{\sin\big( 2n(\theta-\theta_0) \big)}{\sin(\theta-\theta_0)}
  - \frac{\sin\big( 2n(\theta-\theta_0) \big)}{\theta - \theta_0} \,d\theta \bigg|
 \,\le\, \frac{1}{n} \Big( \max_{\theta\in I} h(\theta-\theta_0) \,-\, \min_{\theta\in I} h(\theta-\theta_0) \Big)
\]
for every $I \in \cI_0$ (note that $\theta-\theta_0\in[-\pi/2,\pi/2]$ for $\theta\in I\in\cI_0$).
Summing over intervals $I \in \cI_0$, and noting that we again have an alternating sum, we obtain the bound
\[
\sum_{I \in \cI_0} \bigg| \int_I \frac{\sin\big( 2n(\theta-\theta_0) \big)}{\sin(\theta-\theta_0)} \,
- \, \frac{\sin\big( 2n(\theta-\theta_0) \big)}{\theta - \theta_0} \,d\theta \bigg|
 \le \frac{h(\pi/2) - h(-\pi/2)}{n} = \frac{O(1)}{n},
\]
as claimed.
\end{proof}

It remains to bound $\int_I \frac{\sin(2n(\theta-\theta_0))}{\theta-\theta_0}\,d\theta$ for each $I \in \cI$.
When $d(\theta_0,I) < \pi/n$ we will apply Lemma~\ref{lem:si} to bound this integral. However, in order to deal
with the intervals that are far from $\theta_0$ we will need the following stronger bound. Let $\cJ(\theta_0) := \big\{ I \in \cI_0 : d(\theta_0,I) \ge\pi/n \big\}$.

\medskip
\noindent \textbf{Claim 3:}
$\ds\sum_{I \in \cJ(\theta_0)} \bigg| \ds\int_I \frac{\sin\big( 2n(\theta-\theta_0) \big)}{\theta-\theta_0} \,d\theta \bigg|
 \, \le \, \frac{2}{\pi}$.

\begin{proof}[Proof of Claim~3]
Once again we apply Lemma~\ref{lem:periodic}, this time with $h(\theta)=1/\theta$. We obtain
\[
 \bigg| \int_I \frac{\sin\big( 2n(\theta-\theta_0) \big)}{\theta-\theta_0} \,d\theta \bigg|
 \,\le\, \frac{1}{n} \bigg( \max_{\theta\in I} \frac{1}{\theta-\theta_0}
 \,-\, \min_{\theta\in I} \frac{1}{\theta - \theta_0} \bigg)
\]
for every $I \in \cI_0$ with $\theta_0 \notin I$. Summing over intervals in $\cJ(\theta_0)$,
and noting that we obtain two alternating sums (one on either side of~$\theta_0$), we obtain
\[
 \sum_{I \in \cJ(\theta_0)} \bigg| \ds\int_I \frac{\sin\big( 2n(\theta-\theta_0) \big)}{\theta-\theta_0} \,d\theta \bigg|
 \,\le\, \frac{2}{n} \cdot \frac{1}{\pi/n} = \frac{2}{\pi}
\]
as claimed.
\end{proof}

Note that $2/\pi + 1/(50\pi) + O(1/n) \le 2/3$ if $n$ is sufficiently large, and suppose
first that $\theta_0\in I$ for some $I \in \cI_0$. Then $d(\theta_0,I') \ge \pi/n$ for all $I \ne I' \in \cI_0$,
by Definition~\ref{def:wellseparated}. It follows, by~\eqref{eq:stilde:twostep:again}, Claims 1, 2 and~3,
and Lemma~\ref{lem:si}, that
\[
 \frac{2}{3} = \frac{4}{3} - \frac{2}{3} \le \big|\tg_\alpha(\theta_0)\big| \le 4 + \frac{2}{3} < 5,
\]
as required. On the other hand, if $\theta_0 \notin \bigcup_{I \in \cI_0} I$ then there are at most
two intervals $I \in \cI_0$ such that $d(\theta_0,I)<\pi/n$. Therefore, by~\eqref{eq:stilde:twostep:again},
Claims 1, 2 and~3, and Lemma~\ref{lem:si}, we have
\[
 |\tg_\alpha(\theta_0)| \le 2 \cdot 2 + \frac{2}{3} < 5.
\]
Since $\hat s_\alpha(\theta) = K\sqrt{n}\,\tg_\alpha(\theta)$, this completes the proof of the lemma.
\end{proof}

\begin{remark}
 We note that it is important that the lengths of the intervals $I\in\cI$ are multiples of $\frac{\pi}{n}$.
 Without this assumption it is possible that the error term from the distant intervals $I\in\cI_0$ in Claim~3
 could be unbounded. Indeed, the reason it does not stems ultimately from the cancelation
 in the integrals provided by Lemma~\ref{lem:periodic}.
\end{remark}

Theorem~\ref{thm:sine} is an almost immediate consequence of
Lemmas~\ref{lem:choosing:epsilon} and~\ref{lem:bounding:s:hat}.

\begin{proof}[Proof of Theorem~\ref{thm:sine}]
Let $\cI$ be a suitable and well-separated collection of disjoint intervals in $\R/2\pi\Z$.
By Lemma~\ref{lem:choosing:epsilon}, there exists a colouring $\eps \colon S_o \to \{-1,1\}$
such that, if $\alpha$ is the function given by Lemma~\ref{lem:choosing:alpha}, then
\[
 |s_o(\theta) - \hat s_\alpha(\theta)| \le 72 \sqrt{n}
\]
for every $\theta \in \R$. Now observe that, by Lemma~\ref{lem:bounding:s:hat}, we have
\[
 |s_o(\theta)| \ge |\hat s_\alpha(\theta)| - |s_o(\theta) - \hat s_\alpha(\theta)|
 \ge \bigg( \frac{2K}{3} - 72 \bigg) \sqrt{n} > 10 \sqrt{n}
\]
for all $\theta \in \bigcup_{I\in\cI}I$, and
\[
 |s_o(\theta)| \le |\hat s_\alpha(\theta)| + |s_o(\theta) - \hat s_\alpha(\theta)|
 \le \big( 5K + 72 \big) \sqrt{n} \le 2^{10} \sqrt{n}
\]
for all $\theta\in\R$, as required.
\end{proof}

Finally, let us put together the pieces and prove Theorem~\ref{thm:main:constants}.

\begin{proof}[Proof of Theorem~\ref{thm:main:constants}]
Let $c(\theta)$ be the cosine polynomial, and $\cI$ be the suitable and well-separated
collection of disjoint intervals in $\R/2\pi\Z$, given by Theorem~\ref{thm:cosine}.
Now, given $\cI$, let $s_o(\theta)$ be the sine polynomial given by Theorem~\ref{thm:sine},
and let $s_e(\theta)$ be the sine polynomial defined in~\eqref{def:se}. We claim that the polynomial
\[
 P( e^{i\theta} ) := \big( 1 + 2c(\theta) \big) + 2i \big( s_e(\theta) + s_o(\theta) \big)
\]
has the properties required by the theorem.

To prove the claim, we should first observe that $P(z) = \sum_{k = -2n}^{2n} \eps_k z^k$
with $\eps_k \in \{-1,1\}$ for every $k \in [-2n,2n]$. Indeed the supports
of $c(\theta)$, $s_e(\theta)$, and $s_o(\theta)$ are disjoint and cover the powers
$z^k$ with $k \in \{-2n,\ldots,-2n\} \setminus \{0\}$, and the constant $1$ provides the term corresponding to $k = 0$. Now, observe that
\[
 |P(e^{i\theta})|^2 \le \big(2|c(\theta)|+1\big)^2 + 4|s_e(\theta)+s_o(\theta)|^2
 \le \big(2\sqrt{n}+1\big)^2 + 4\big(2^{10}+6\big)^2 n
 \le \big(2^{12}\sqrt{n}\big)^2,
\]
for every $\theta \in \R$, since $|c(\theta)| \le \sqrt{n}$ and $|s_e(\theta)| \le 6\sqrt{n}$,
by Theorem~\ref{thm:cosine} and Lemma~\ref{lem:even_ubounds}, and $|s_o(\theta)| \le 2^{10} \sqrt{n}$,
by Theorem~\ref{thm:sine}. Next, observe that if $\theta \notin \bigcup_{I\in\cI} I$ then
\[
 |P(e^{i\theta})| \ge \big| \Re\big(P(e^{i\theta})\big) \big|
 \ge 2|c(\theta)|-1 \ge \delta \sqrt{n}
\]
for all sufficiently large~$n$, by Theorem~\ref{thm:cosine}.
Finally, if $\theta \in \bigcup_{I \in \cI} I$, then
\[
 |P(e^{i\theta})| \ge \big| \Im\big(P(e^{i\theta}) \big)\big|
 \ge 2\big(|s_o(\theta)|-|s_e(\theta)|\big)
 \ge 2\big(10\sqrt{n}-6\sqrt{n}\big)=8\sqrt n,
\]
by Theorem~\ref{thm:sine}. Hence $|P(z)|\ge\delta\sqrt{n}$ for all $z\in\C$ with $|z|=1$, as required.
\end{proof}

\section*{Acknowledgements}

Much of this research was carried out during a one-month visit by the authors to IMT Lucca. We are grateful to IMT (and especially to Prof.~Guido Caldarelli) for providing a wonderful working environment.


\begin{thebibliography}{99}

\bibitem{PB} P.~Balister,
Bounds on Rudin--Shapiro polynomials of arbitrary degree,
\emph{in preparation.}

\bibitem{B81} J.~Beck,
Roth's estimate of the discrepancy of integer sequences is nearly sharp,
\emph{Combinatorica}, {\bf 1} (1981), 319--325.

\bibitem{B91} J.~Beck,
Flat polynomials on the unit circle -- note on a problem of Littlewood,
\emph{Bull.~London Math.~Soc.}, {\bf 23} (1991), 269--277.

\bibitem{BN} E.~Beller and D.J.~Newman,
The Minimum Modulus of Polynomials,
\emph{Proc.~Amer.~Math.~Soc.}, {\bf 45} (1974), 463--465.

\bibitem{BP32} A.~Bloch and G.~P\'olya,
On the roots of certain algebraic equations,
\emph{Proc.~Lond.~Math.~Soc.}, {\bf 33} (1932), 102--114.

\bibitem{BB} E.~Bombieri and J.~Bourgain,
On Kahane's ultraflat polynomials,
\emph{J.~Eur.~Math.~Soc.}, {\bf 11} (2009), 627--703.

\bibitem{B86} J.~Bourgain, 
Sur le minimum d'une somme de cosinus, 
\emph{Acta Arith.}, {\bf 45} (1986), 381--389.

\bibitem{B02} P.~Borwein,
Computational Excursions in Analysis and Number Theory, Springer-Verlag, New York, 2002. 

\bibitem{B73} J.~Brillhart,
On the Rudin--Shapiro polynomials,
\emph{Duke Math.~J.}, {\bf 40} (1973), 335--353.

\bibitem{BC} J.~Brillhart and L.~Carlitz,
Note on the Shapiro polynomials,
\emph{Proc.~Amer.~Math.~Soc.}, {\bf 25} (1970), 114--119.


\bibitem{BLM} J.~Brillhart, J.S.~Lomont and P.~Morton,
Cyclotomic properties of the Rudin--Shapiro polynomials,
\emph{J.~Rein.~Angew.~Math.}, {\bf 288} (1976), 37--65.

\bibitem{B77} J.S.~Byrnes,
On polynomials with coefficients of modulus one,
\emph{Bull.~London Math.~Soc.}, {\bf 9} (1977), 171--176.

\bibitem{CEF} F.W.~Carroll, D.~Eustice and T.~Figiel,
The minimum modulus of polynomials with coefficients of modulus one,
\emph{J.~London Math.~Soc.}, {\bf 16} (1977), 76--82.

\bibitem{C65} S.~Chowla, 
Some applications of a method of A.~Selberg, 
\emph{J.~Reine Angew.~Math.}, \textbf{217} (1965), 128--132.


\bibitem{E57} P.~Erd\H{o}s,
Some unsolved problems,
\emph{Michigan Math.~J.}, {\bf 4} (1957), 291--300.

\bibitem{EO56} P.~Erd\H{o}s and A.C.~Offord,
On the number of real roots of a random algebraic equation,
\emph{Proc.~London Math.~Soc.}, {\bf 6} (1956), 139--160.

\bibitem{G49} M.J.E.~Golay,
Multislit spectrometry,
\emph{J.~Opt.~Soc.~Am.}, {\bf 39} (1949), 437--444.

\bibitem{HL16} G.H.~Hardy and J.E.~Littlewood,
Some problems of Diophantine approximation: a remarkable trigonometric series,
\emph{Proc.~Nat.~Acad.~Sci.}, {\bf 2} (1916), 583--586.

\bibitem{HL48} G.H.~Hardy and J.E.~Littlewood, 
A new proof of a theorem on rearrangements, 
\emph{J.~London Math.~Soc.}, \textbf{23} (1948), 163--168.

\bibitem{H91} G.~Howell,
Derivative error bounds for Lagrange interpolation:
an extension of Cauchy's bound for the error of Lagrange interpolation,
\emph{J.~Approx.~Theory}, {\bf 67} (1991), 164--173.

\bibitem{Ka80} J.-P.~Kahane,
Sur les polyn\^omes \`a coefficients unimodulaires,
\emph{Bull.~London Math.~Soc.}, {\bf 12} (1980), 321--342.

\bibitem{K81} S.V.~Konyagin, On a problem of Littlewood, 
\emph{Izv.~Akad.~Nauk SSSR Ser.~Mat.}, {\bf 45} (1981), 243--265. 

\enlargethispage*{\baselineskip}
\enlargethispage*{\baselineskip}
\thispagestyle{empty}

\bibitem{Ko80} T.W.~K\"orner,
On a polynomial of Byrnes,
\emph{Bull.~London Math.~Soc.}, {\bf 12} (1980), 219--224.

\bibitem{L61} J.E.~Littlewood,
On the mean values of certain trigonometric polynomials,
\emph{J.~London Math.~Soc.} {\bf 36} (1961), 307--334.

\bibitem{L62a} J.E.~Littlewood,
On the mean values of certain trigonometric polynomials~II,
\emph{Illinois J.~Math.}, {\bf 6} (1962), 1--39.

\bibitem{L62} J.E.~Littlewood,
On the real roots of real trigonometrical polynomials, In:~Studies in Mathematical Analysis and Related Topics:~Essays in Honor of George P\'olya (G.~Szeg\H{o}, ed.), pp.~219--226, Stanford Univ.~Press, Stanford, Calif., 1962.

\bibitem{L64} J.E.~Littlewood,
On the real roots of real trigonometrical polynomials, II,
\emph{\em J.~London Math.~Soc.}, {\bf 39} (1964), 511--532.

\bibitem{L66a} J.E.~Littlewood,
The real zeros and value distributions of real trigonometrical polynomials,
\emph{J.~London Math.~Soc.}, {\bf 41} (1966), 336--342.

\bibitem{L66} J.E.~Littlewood,
On polynomials $\sum^n \pm z^m$, $\sum^n e^{\alpha_m i} z^m$, $z=e^{\theta i}$,
\emph{J.~London Math.~Soc.}, {\bf 41} (1966), 367--376.

\bibitem{L68} J.E.~Littlewood,
Some Problems in Real and Complex Analysis,
D.C.~Heath and Co., Raytheon Education Co., Lexington, Mass, 1968.

\bibitem{LO38} J.E.~Littlewood and A.C.~Offord,
On the number of real roots of a random algebraic equation,
\emph{J.~London Math.~Soc.}, {\bf 13} (1938),  288--295.

\bibitem{LO45} J.E.~Littlewood and A.C.~Offord,
On the distribution of zeros and $a$-values of a random integral function~I,
\emph{J.~London Math.~Soc.}, {\bf 20} (1945), 130--136.

\bibitem{LO48} J.E.~Littlewood and A.C.~Offord,
On the distribution of zeros and $a$-values of a random integral function~II,
\emph{Ann.~Math.}, {\bf 49} (1948), 885--952; errata {\bf 50} (1949), 990--991.

\bibitem{LM} S.~Lovett and R.~Meka,
Constructive Discrepancy Minimization by Walking on The Edges,
\emph{SIAM J.~Computing}, {\bf 44} (2015), 1573--1582.


\bibitem{MPS} O.C.~McGehee, L.~Pigno and B.~Smith, 
Hardy's inequality and the $L_1$ norm of exponential sums, 
\emph{Ann.~Math}, \textbf{113} (1981), 613--618.

\bibitem{M17} H.L.~Montgomery,
Littlewood polynomials,
In:~Analytic Number Theory, Modular Forms and $q$-hypergeometric Series (G.~Andrews and F.~Garvan, eds.), pp.~533--553, Springer, Cham, 2017.

\bibitem{O18} A.~Odlyzko,
Search for ultraflat polynomials with plus and minus one coefficients,
In:~Connections in Discrete Mathematics:~A Celebration of the Work of Ron Graham
(S.~Butler, J.~Cooper and G.~Hurlbert, eds.), pp.~39--55, Cambridge Univ.~Press, 2018.

\bibitem{R17} B.~Rodgers,
On the distribution of Rudin--Shapiro polynomials and lacunary walks on $SU(2)$,
\emph{Adv.~Math.}, {\bf 320} (2017), 993--1008.

\bibitem{R59} W.~Rudin,
Some theorems on Fourier coefficients,
\emph{Proc.~Amer.~Math.~Soc.}, {\bf 10} (1959), 855--859.

\bibitem{R04} I.Z.~Ruzsa, Negative values of cosine sums, 
\emph{Acta Arith.}, {\bf 111} (2004), 179--186.

\bibitem{SZ54} R.~Salem and A.~Zygmund,
Some properties of trigonometric series whose terms have random signs,
\emph{Acta Math.} {\bf 91} (1954), 245--301.

\bibitem{S41} A.C.~Schaeffer, 
Inequalities of A.~Markoff and S.~Bernstein for polynomials and related functions, 
\emph{Bull.~Amer.~Math.~Soc.}, {\bf 47} (1941), 565--579.

\enlargethispage*{\baselineskip}
\enlargethispage*{\baselineskip}
\thispagestyle{empty}

\bibitem{S52} H.S.~Shapiro,
Extremal problems for polynomials,
Thesis for S.M.~Degree, 1952, 102 pp.

\bibitem{S85} J.~Spencer,
Six standard deviations suffice,
\emph{Trans.~Amer.~Math.~Soc.}, {\bf 289} (1985), 679--706.

\end{thebibliography}
\end{document}